\documentclass{article}
\usepackage{amssymb}

\usepackage{amsmath}

\advance\textheight by \topskip
\begin{document}

\title{A Theorem of Legendre in $I\Delta_0+\Omega_1$}

\author{Paola D'Aquino and Michele Bovenzi}

\maketitle

\newtheorem{fact}{Fact}[section]

\newtheorem{theor}[fact]{Theorem}

\newtheorem{lem}[fact]{Lemma}

\newtheorem{cor}[fact]{Corollary}

\newtheorem{defn}[fact]{Definition}

\newtheorem{remark}[fact]{Remark}

\newtheorem{prop}[fact]{Proposition}

\newtheorem{rem}[fact]{Remark}

\begin{abstract}

We prove a classical theorem due to Legendre, about the existence of non trivial solutions of quadratic diophantine equations of the form $ax^2+by^2+cz^2=0$, in the weak fragment of Peano Arithmetic $I\Delta_0+\Omega_1$.

\end{abstract}

\section{Introduction}

\indent 

We work with the language $\cal L$ of arithmetic containing the symbols $\{ 0,1,+,\cdot,\leq \}$ and we focus on the theory of {\it bounded induction} $I\Delta_0$, the fragment of Peano Arithmetic ($PA$) where induction is restricted only to bounded formulas ($\Delta_0$-formulas).

It is well known that $I\Delta_0$ does not prove the totality of the exponential function (see \cite{Parikh}). This is the main limitation for reproducing in $I\Delta_0$ many classical results of elementary number theory that rely on functions of exponential growth. For example, it is still unknown if $I\Delta_0$ proves the existence of arbitrarily large of primes.

A much stronger theory is obtained if we add an axiom, denoted by $exp$, which guarantees the totality of the exponential function. The theory $I\Delta_0+exp$ is strong enough to 
reproduce almost all elementary number theory (as, for example, from \cite{HW}). Woods (\cite{Woods}) proved that in some cases functions of exponential growth can be avoided and replaced by some combinatorial principle, such as the {\it pigeonhole principle}. In fact, Woods proved unboundedness of primes in the theory $I\Delta_0+\Delta_0PHP$, where $\Delta_0PHP$ denotes an axiom stating a $\Delta_0$-version of the pigeonhole principle (namely that there exist no injective $\Delta_0$-definable function from $z+1$ into $z$).

Later, Paris, Wilkie and Woods (\cite{Paris-Wilkie-Woods}) improved this result by showing that  a {\it weaker} version of the PHP, denoted by $\Delta_0$-WPHP, is sufficient in order to prove unboundedness of primes. Moreover, they showed that $\Delta_0$-WPHP is provable in the theory $I\Delta_0+\Omega_1$, where $\Omega_1$ is the axiom $\forall x\forall y\exists z\ \left( x^{log_{_{2}}y}=z\right)$. Hence, $I\Delta_0+\Omega_1$ proves cofinality of primes.

The theories $I\Delta_0$ and  $I\Delta_0+\Omega_1$ have been widely studied also for their connections with complexity theory (see \cite{Wilkie}). Many open problems in $I\Delta_0$ or $I\Delta_0+\Omega_1$ have complexity-theoretic counterparts. For example, it is still an open problem if $I\Delta_0$ proves the {\it MRDP}-theorem (from Matijasevic, Robinson, Davis and Putnam). This theorem asserts that every recursively enumerable set is existentially definable.
A formalization of it in the language $\cal L$ is:\\
for any $\Sigma_1$-formula does it exist a polynomials $p({\bar x},{\bar y})$ over $\mathbb Z$ such that
$$ I\Delta_0\vdash\ \forall{\bar x}\ \left(\theta({\bar x})\leftrightarrow\exists{\bar y}\ p({\bar x},{\bar y})=0\right)?$$
Wilkie (see \cite{Wilkie}) observed that a positive answer to this problem in $I\Delta_0$ would give a positive solution to the well known open problem in complexity theory if $NP\overset{?}{=}coNP$. It is also unknown if $I\Delta_0+\Omega_1$ proves $MRDP$-theorem and, again, a positive answer would give $NP=coNP$.
 
On the other hand, many other classical number-theoretical properties have been proved in $I\Delta_0+\Omega_1$, such as Lagrange's four squares theorem (see \cite{Ber-Intr}) and basic results about residue fields (see \cite{D-M}). In this paper we show that $I\Delta_0+\Omega_1$ proves a classical theorem due to Legendre about quadratic diophantine equations.

\subsection{Bounded induction and $I\Delta_0+\Omega_1$}

In this section we recall some basic properties of the theories $I\Delta_0$ and $I\Delta_0+\Omega_1$ that have been used for the main proofs in Chapter \ref{chapterTheorem}.\\

We recall that in $I\Delta_0$ induction is allowed only on formulas in which all quantifiers are bounded by terms of the language. Notice that in the language of arithmetic, terms are actually polynomials. 

When working in $I\Delta_0$ care must be taken in expressing properties via $\Delta_0$-formulas, since these are the only formulas we can induct on. For example the very basic statement about divisibility is expressed in a $\Delta_0$-way as follows
\begin{equation}\label{D0divides}
x\mbox{ divides } y:\ \delta(x,y)=\exists z\leq y\ (xz=y).
\end{equation}
By $\Delta_0$-induction it is easily proved that any two elements in a model of $I\Delta_0$ have a greatest common divisor, and this can be expressed by the Bezout identity, and thus all models of $I\Delta_0$ are Bezout rings. Moreover, $\Delta_0$-induction can be used to prove that any non empty set which is $\Delta_0$-definable has a minimum element.

As we already recalled, it is still unknown if any models of $I\Delta_0$ has cofinally many primes. This does not affect the factorization in powers of primes of any element. In order to see this we use some lemmas.

\begin{lem}\label{boundedset}
Let $A$ be a bounded and $\Delta_0$-definable subset of ${\cal M}\models I\Delta_0$. Then $A$ has a maximum element.
\end{lem}

{\it{Proof.}} Let $\phi(x)$ be the $\Delta_0$ formula defining the set $A$, and let $\alpha\in\cal M$ be an upper bound for $A$.\\
The set $X=\{ y\leq\alpha\ :\ \exists t\leq\alpha\ (\phi(t)\land y<t)\}$ is clearly $\Delta_0$-definable, and so is the set $X^c={\cal M}\setminus X=\{ y\ :\ y>\alpha\}\cup\{ y\ :\ y\leq\alpha\land\forall t\leq\alpha ( \phi(t)\rightarrow y>t)\}$. Let $x_0=min(X^c)$; then $x_0-1\in X$, hence there is $t\in A$ such that $x_0-1\leq t$. But now necessarily $x_0-1=t$, so $x_0-1=max(A)$. $_\Box$

\begin{lem}\label{lcm}
Let $A$ be a bounded, $\Delta_0$-definable subset of ${\cal M}\models I\Delta_0$. If there is a non-zero $m\in {\cal M}$ divisible by all $a\in A$, then there is a non-zero $\mu\in {\cal M}$ which is minimal with respect to this property (i.e. if $x\in{\cal M}$ is divisible by all elements of $A$, then $\mu$ divides $x$).
\end{lem}

{\it{Proof.}} Let $\phi(x)$ be the $\Delta_0$-formula defining the set $A$.\\
The set $D=\{ b\ :\ b\not= 0\land \forall a\leq m\ (\phi(a)\rightarrow \delta(a,b))\}$ is $\Delta_0$-definable and nonempty since $m\in D$ (here $\delta(a,b)$ is the $\Delta_0$-formula for divisibility (\ref{D0divides})).

Let $\mu$ be the minimum of $D$. Now let $x\in{\cal M}$ be divisible by all elements of $A$ and not by $\mu$, then we have $x=\mu q+r$ for $q,r\in{\cal M}$, with $0\leq r< \mu$. Hence $r=x-\mu q$ is divisible by all elements of $A$, and since $\mu=min(D)$, it must be $r=0$ and so $\mu$ divides $x$. $_\Box$

The minimal element $\mu$ divisible by all elements of $A$ will be called the {\it least common multiple of $A$} and it will be denoted by $lcm(A)$.

We also prove, by an easy $\Delta_0$-induction, that every element is divisible by a prime.

\begin{lem}
$I\Delta_0\vdash \forall x>1\ \exists p\leq x\ (Pr(p)\land \delta(p,x))$.
\end{lem}

{\it{Proof.}} Here $Pr(x)$ is the $\Delta_0$-statement for primes
\begin{equation}\label{D0prime}
Pr(x)=\forall y\leq x\ (\delta(y,x)\rightarrow (y=1\lor y=x)).
\end{equation}

Let $\psi(x)=\forall y\leq x\ \exists p\leq y\ (Pr(p)\land \delta(p,y))$; clearly $I\Delta_0\vdash \psi(2)$. Suppose $I\Delta_0\vdash \psi(x)$ and consider $x+1$. If $x+1$ is a prime, then $I\Delta_0\vdash \psi(x+1)$. If $x+1$ is not a prime, then there is a $1< y\leq x$ that divides $x+1$. But from $I\Delta_0\vdash \psi(x)$ we deduce that there is a prime dividing $y$, and thus dividing $x+1$, hence $I\Delta_0\vdash \psi(x+1)$. $_\Box$

We can express that an element is a power of a prime $p$ with the $\Delta_0$-formula
\begin{equation}\label{powerprimes}
Pow_p(x) = Pr(p)\land\forall y\leq x\ (\delta(y,x)\rightarrow\delta(p,y)).
\end{equation}
This allows us to identify the greatest power of a prime that divides a given element.

\begin{lem}
$I\Delta_0\vdash \forall x>1\ \forall p\ (Pr(p)\rightarrow \exists!y\leq x\ (Pow_p(y)\land\delta(y,x)\land\lnot\delta(py,x))$
\end{lem}

{\it{Proof.}}  Let ${\cal M}\models I\Delta_0$, $x,p\in{\cal M}$, with $x>1$ and $p$ a prime.\\
The set $A_p^x=\{ y\ :\ Pow_p(y)\land \delta(y,x)\}$ of all powers of $p$ dividing $x$ is clearly $\Delta_0$-definable and bounded by $x$, hence by Lemma \ref{boundedset} it has a maximum element. $_\Box$
\medskip

We can now express that "$y$ is the greatest power of $p$ that divides $x$" with the $\Delta_0$-formula
\begin{equation}
MPow_p(x,y)=Pow_p(y)\land\delta(y,x)\land\ \forall z\leq x\ (Pow_p(z)\land\delta(z,x))\rightarrow \delta(z,y)
\end{equation}

For any $x\in{\cal M}\models I\Delta_0$ we can consider the $\Delta_0$-definable set
\begin{equation}
A_x=\{ y\ :\ \exists p\leq x\ (Pr(p)\land MPow_p(x,y))\}
\end{equation}
of all maximum powers of primes dividing $x$. Since $x\in{\cal M}$ is clearly divisible by all elements of $A_x$, by Lemma \ref{lcm} there is the smallest $\mu$ divisible by all elements of $A_x$, which is trivially shown to coincide with $x$. Hence we have the property of factorization in powers of primes for every element $x$ of a model of $I\Delta_0$ as the $lcm(A_x)$.

What will also be used later is the following $I\Delta_0$-version of the {\it Chinese Reminder Theorem (CRT)} (see \cite{D'Aquino}).

\begin{theor}\label{D0CRT}
Let $A$ be a bounded, $\Delta_0$-definable subset of ${\cal M}\models I\Delta_0$. Let $f,r:A\longrightarrow{\cal M}$ be a $\Delta_0$-definable functions such that $(f(a),f(b))=1$ for every $a,b\in{\cal M}$ and $r(a)<f(a)$ for every $a\in{\cal M}$. Suppose there is $w\in{\cal M}$ which is divisible by all elements of $f(A)$.\\
Then there is $u<\prod_{a\in A}f(a)$ such that $u\equiv r(a)(mod\ f(a))$ for every $a\in A$.
\end{theor}

\begin{remark}
\end{remark}\vspace{-1em}Notice that we can express in models of $I\Delta_0$ the congruence "$x$ is equivalent to $y$ modulo $z$" via the $\Delta_0$-formula:
\begin{equation}\label{D0equiv}
x\equiv y (mod\ z)\ :\ \exists k\leq x\ (x=y+kz).
\end{equation}
Now we remark that the previous theorem is a generalization of the classic CRT where $A=\{m_1,\dots,m_k\}\subseteq {\mathbb Z}$ is a finite set of pairwise relatively prime moduli, $r_i<m_i$ for all $i=1,\dots,k$ and we are looking for integer solutions of the set of congruences
$$\begin{cases}
x\equiv r_1 (mod\ m_1) \\ \ \ \ \vdots \\ x\equiv r_k (mod\ m_k)
\end{cases}.$$
In a model $\mathcal M$ of $I\Delta_0$, the $I\Delta_0$-CRT allows us to extend such property to any bounded, $\Delta_0$-definable subset $A$ of pairwise relatively prime moduli.\\

\medskip

Even though the theory $I\Delta_0$ is strong enough to prove factorization of elements as products of primes, there are many other classical number-theoretical results whose provability in $I\Delta_0$ are still open problems. The main obstacle in obtaining results like unboundedness of primes is, as we mentioned before, the lack of functions of exponential growth rate, such as factorials, since they are not provably total in models of $I\Delta_0$.

As mentioned before, Woods (see \cite{Woods}) proved that in some cases such functions can be avoided and replaced by combinatorial principles such as the {\it Pigeonhole Principle}. He showed that if we add to the theory $I\Delta_0$ the following $\Delta_0$-version of the pigeonhole principle ($\Delta_0$-PHP)
$$
\forall x<z\ \exists y<z\ \theta(x,y)\rightarrow \exists x_1<z+1\ \exists x_2< z+1\ (x_1\not= x_2\land \theta(x_1,y)\land\theta(x_2)),
$$
where $\theta(x,y)$ runs through all $\Delta_0$-formulas, the resulting theory $I\Delta_0+\Delta_0$-PHP is strong enough to prove the existence of arbitrarily large primes.

\begin{theor}[Woods, \cite{Woods}]
$$I\Delta_0+\Delta_0\mbox{-PHP}\vdash \forall x\exists p\ (Pr(p)\land p>x).$$
\end{theor}

Actually a weaker version of the PHP turned out to be sufficient to prove unboundedness of primes in $I\Delta_0$. Paris, Wilkie and Woods in \cite{Paris-Wilkie-Woods} used such a principle, denoted by $\Delta_0$-WPHP, asserting that there is no injective $\Delta_0$-function from $(1+\varepsilon)z$ into $z$, for every rational number $\varepsilon$ such that the integer part of $(1+\varepsilon)z$ is greater than $z$. This latter result is hence as follows.

\begin{theor}[Paris, Wilkie and Woods, \cite{Paris-Wilkie-Woods}]
$$I\Delta_0+\Delta_0\mbox{-WPHP}\vdash \forall x\exists p\ (Pr(p)\land p>x).$$
\end{theor}

In \cite{Paris-Wilkie-Woods} it is also shown that the $\Delta_0$-WPHP is provable in the theory $I\Delta_0+\Omega_1$, where $\Omega_1$ is the axiom
\begin{equation}\label{omega1}
\forall x\forall y\exists z\ \left( x^{log_{_{2}}y}=z\right).
\end{equation}

Notice that the quantity $log_{_{2}}y$ has a $\Delta_0$-meaning, as the following lemma guarantees.

\begin{lem}\label{D0log}
Let ${\cal M}\models I\Delta_0$ and let $a,m\in{\cal M}$, with $m>0$, $a>1$. Then there is a unique $l_0\in{\cal M}$ such that $a^{l_0}\leq m < a^{l_0+1}$.
\end{lem}

{\it{Proof.}}
The set $A=\{ l\in{\cal M}\ :\ \exists b\leq m\ E_0(a,l,b)\}$ is clearly $\Delta_0$-definable in $\cal M$, and it is bounded and non-empty. Hence, by Lemma \ref{boundedset}, it has a maximum element $l_0$, so it is $a^{l_0}\leq m < a^{l_0+1}$. $_\Box$

We can state the following.

\begin{theor}[Paris, Wilkie and Woods, \cite{Paris-Wilkie-Woods}]
$I\Delta_0+\Omega_1\vdash \forall x\exists p\ (Pr(p)\land p>x).$
\end{theor}

\section{Legendre's Theorem}\label{chapterTheorem}

\indent 

Our proof of Legendre's theorem in $I\Delta_0+\Omega_1$ follows the lines of the corresponding theorem given in \cite{IR}. Our contribution has involved a careful analysis of the objects used in the proof. In particular, we have ensured that all the properties and tools involved in the proof are $\Delta_0$-definable and valid in our theory. We have finally obtained estimates on the growth rate of the solution of the considered equations, of polynomial size in  $x^{log_{_2} y}$, hence proving the main theorem in $I\Delta_0+\Omega_1$ by $\Delta_0$-induction.

\subsection{The theorem and its equivalent}

\indent

Legendre's theorem gives a necessary and sufficient condition for the existence of non trivial solution for certain quadratic diophantine equations. The equations considered are of the form
\begin{equation}\label{eq.legendre}
ax^2+by^2+cz^2=0,
\end{equation}
with $a,b,c\in \mathbb Z\setminus\{0\}$, {\it square free}, relatively prime and, clearly, not all of the same sign. For a non trivial solution we mean a solution $(x_0,y_0,z_0)\in \mathbb Z$ different from $(0,0,0)$.

An integer $a$ is {\it square free} if no square divides $a$. This property is $\Delta_0$-definable via the formula:
$$
\sigma(x)=\forall y\leq x\ (Pr(y)\rightarrow \lnot\delta(y^2,x)),
$$
where $\delta(y,x)=${\it "y divides x"} and $Pr(y)=${\it "y is a prime"} are expressed by the $\Delta_0$-formulas:
\begin{equation}
\delta(y,x)= \exists z\leq x\ (yz=x)
\end{equation}
\begin{equation}
Pr(y)= \forall x\leq y\ (\delta(x,y)\rightarrow (x=1\lor y=x)).
\end{equation}

We will consider {\it congruences} in $I\Delta_0$: as for the classical definition on integers we say that $a$ is congruent to $b$ modulo $c$, with $a,b,c\in{\cal M}\models I\Delta_0$, $c\not=0$, if $c$ divides $b-a$, and this can be expressed by the $\Delta_0$-formula
\begin{equation}\label{D0cong}
\gamma(a,b,c)=\exists m\leq a\ (a=mc+b).
\end{equation}
We will also denote this fact as $a\equiv b(mod\ c)$. Since euclidean division is valid in $I\Delta_0$, given any $a,c\in{\cal M}$, $c\not=0$, we can find $q,r\in{\cal M}$ such that $a=qc+r$, with $0\leq r<c$, and so we have $a\equiv r\equiv r-c(mod\ c)$, with either $r$ or $r-c\leq c/2$. Therefore we can always consider that $a\equiv b(mod\ c)$ with $|b|\leq c/2$.

Models of $I\Delta_0$ are Bezout rings, so if $(a,b)=1$ then there are $h, k$ such that $ah+bk=1$, so $ah\equiv 1(mod\ b)$, and we can identify $a^{-1}$ with $h$. This is a $\Delta_0$-property since it can be expressed by the $\Delta_0$-formula
\begin{equation}\label{D0invert}
Inv(a,b)=\exists w\leq b\ \gamma(aw,1,b),
\end{equation}
where $\gamma$ is the $\Delta_0$-formula (\ref{D0cong}).

We will denote the fact that {\it a is a square modulo b} by $a{\cal R}b$. This is $\Delta_0$-definable via:
$$
\rho(a,b)=\exists x\leq b/2 \land \exists m\leq a (a=x^2+mb).
$$

The statement of Legendre's theorem in a model ${\cal M}$ of $I\Delta_0+\Omega_1$ is as follows.

\begin{theor}[Legendre]\label{Legendre}
Let $a,b,c\in {\cal M}\setminus\{0\}$ be non-zero, not all of the same sign, square-free and pairwise relatively prime.
Then the equation
\begin{equation}\label{eq.legendre.2}
ax^2+by^2+cz^2=0
\end{equation}
has a non trivial solution in ${\cal M}$ if and only if
\begin{description}
 \item{\rm (\it Leg.1)}\ $-ab{\cal R}c$,
 \item{\rm (\it Leg.2)}\ $-bc{\cal R}a$,
 \item{\rm (\it Leg.3)}\ $-ac{\cal R}b$.
\end{description}
\end{theor}

\medskip

\begin{remark}
\end{remark}\vspace{-1em}
{\bf 1.} If a non trivial solution $(x_0,y_0,z_0)$ of (\ref{eq.legendre.2}) exists, then we can consider $x_0,y_0,z_0$ to be pairwise relatively prime (and call such a solution {\it primitive}). Indeed, if $p$ is a prime which divides, say, $x_0$ and $y_0$, from $ax_0^2+by_0^2+cz_0^2=0$, and being $c$ square free, it follows that $p$ also divides $z_0$, so we can factor out $p$ and consider the solution $(x_0/p,y_0/p,z_0/p)$.

\noindent
{\bf 2.} The necessary condition of Theorem \ref{Legendre} is proved as follows.

{\it{Proof.}} ($\Longrightarrow$ of \ref{Legendre}) Let $(x_0,y_0,z_0)$ be a primitive solution of (\ref{eq.legendre}) in a model $\cal M$ of $I\Delta_0$. Then
\begin{equation}\label{eq1}
ax_0^2+by_0^2+cz_0^2=0,
\end{equation}
and reducing modulo $a$ we get $-by_0^2\equiv cz_0^2(mod\ a)$.\\
Now if a prime $p$ divides $a$ and $y_0$, equation (\ref{eq1}) implies that $p|cz_0^2$. From $(a,c)=1$ we get that $p\not| c$ and so $p|z_0$. Hence $p$ divides $y_0$ and $z_0$, which is a contradiction.

Hence we have $(a,y_0)=1$, and $y_0$ is invertible modulo $a$. So we can write $-b\equiv cz_0^2(y_0^{-1})^2 (mod\ a)$, (here $y_0^{-1}$ is given by (\ref{D0invert})). Hence we can rewrite $-bc\equiv(cz_0y_0^{-1})^2 (mod\ a)$, i.e. $-bc{\cal R}a$.
In the same way, being the equation (\ref{eq.legendre}) symmetric for $a,b$ and $c$, we also get $ -ab{\cal R}c$ and $-ac{\cal R}b$. $_\Box$\\

\medskip

We will prove Legendre's theorem in $I\Delta_0+\Omega_1$ in the following equivalent form, which we call {\it normal form}. Using the same notation as before, the statement is as follows.

\begin{theor}[Legendre normal form]\label{Legendre.normal}
Let $a,b\in {\cal M}\setminus\{0\}$ be square-free and positive.
Then the equation
\begin{equation}\label{eq.legendre.normal}
ax^2+by^2=z^2
\end{equation}
has a non trivial solution in ${\cal M}$ if and only if
\begin{description}\vspace{-0.5em}
 \item{\rm (\it Norm.1)}\ $ a{\cal R}b$,
 \item{\rm (\it Norm.2)}\ $b{\cal R}a$,
 \item{\rm (\it Norm.3)}\ $-{\frac{ab}{d^2}}{\cal R}d,$ where $d=(a,b)$.
\end{description}
\end{theor}

\begin{remark}
\end{remark}\vspace{-1em}
{\bf 1.} Using the same arguments as before we can assume that, if a non trivial solution $(x_0,y_0,z_0)$ of the equation (\ref{eq.legendre.normal}) exists, it is {\it primitive} (i.e. $x_0,y_0,z_0$ pairwise relatively coprime).

\noindent
{\bf 2.} The necessary condition of the Theorem \ref{Legendre.normal} is proved as follows.

{\it{Proof.} ($\Longrightarrow$ of \ref{Legendre.normal})} Let $(x_0,y_0,z_0)$ be a primitive solution of (\ref{eq.legendre.normal}) in a model $\cal M$ of $I\Delta_0$, i.e.
\begin{equation}\label{eq1.normal}
ax_0^2+by_0^2=z_0^2.
\end{equation}

Since the solution is primitive, it can be easily deduced that $(a,y_0)=1$. So $y_0$ is invertible modulo $a$ in $\cal M$, and we have $by_0^2\equiv z_0^2(mod\ a)$, which implies $b\equiv (z_0y_0^{-1})^2(mod\ a)$, i.e. $b{\cal R}a$.
Similarly, by symmetry of the equation (\ref{eq1.normal}) on the coefficients $a$ and $b$, we can show that $a{\cal R}b$.

In order to obtain condition {\rm (\it Norm.3)} of the theorem, we first observe that $d$ is square-free since $a,b$ are. Moreover, if $a=da'$ and $b=db'$, with $(a',b')=1$, from (\ref{eq1.normal}) we get
\begin{equation}\label{eq2.normal}
d(a'x_0^2+b'y_0^2)=z_0^2,
\end{equation}
i.e. $d|z_0^2$. Since $d$ is square-free $d$ divides $z_0$. Let $z_0=dz_0'$.\\
From (\ref{eq2.normal}) we obtain
\begin{equation}\label{eq4.normal}
\frac{a}{d}x_0^2+\frac{b}{d}y_0^2=d\left(\frac{z_0}{d}\right)^2.
\end{equation}
Hence $-\frac{a}{d}x_0^2=\frac{b}{d}y_0^2-d\left(\frac{z_0}{d}\right)^2\equiv \frac{b}{d}y_0^2(mod\ d)$.

We already noticed that $(x_0,b)=1$, so also $(x_0,d)=1$, hence $x_0$ is invertible modulo $d$. So we have
$$
-\frac{a}{d}\equiv \frac{b}{d}y_0^2(x_0^{-1})^2(mod\ d) \Rightarrow -\frac{a}{d}\frac{b}{d}\equiv \left(\frac{b}{d}y_0x_0^{-1}\right)^2(mod\ d),
$$
that means $-{\frac{ab}{d^2}}{\cal R}d$.$_\Box$\\

\medskip

Now we prove the equivalence between the two statements of Legendre's theorem. We will use the following lemma.

\begin{lem}\label{lemmaCRT}
Let $a,m,n\in{\mathcal M}\models I\Delta_0$, $m,n$ relatively prime. If $a{\cal R}m$ and $a{\cal R}n$, then $a{\cal R}mn$.
\end{lem}

{\it{Proof.}} From the hypothesis there are $\alpha,\beta\in{\cal M}$, $\alpha\leq m/2, \beta\leq n/2$ such that
$$a\equiv \alpha^2(mod\ m) \mbox{ and } b\equiv\beta^2(mod\ n).$$

Consider the system of congruences
\begin{equation}\label{sys}
\begin{cases}
x\equiv \alpha (mod\ m) \\
x\equiv \beta (mod\ n)
\end{cases}.
\end{equation}
Since $(m,n)=1$, the system (\ref{sys}) has a solution $\gamma\in\cal M$ by $\Delta_0$-CRT. Hence, we have
$$
\gamma\equiv \alpha (mod\ m) \mbox{ and } \gamma\equiv \beta (mod\ n),
$$
And so we obtain
$$
a\equiv \gamma^2 (mod\ m) \mbox{ and } a\equiv \gamma^2 (mod\ n),
$$
that means that both $m,n$ divide $a-\gamma^2$. Since $m,n$ are relatively prime, $mn$ divides $a-\gamma^2$, that means $a\equiv \gamma^2(mod\ mn)$, i.e. $a{\cal R}mn$.$_\Box$

\medskip

\begin{theor}\label{equivalence}
Theorems \ref{Legendre} and \ref{Legendre.normal} are equivalent.
\end{theor}

{\it{Proof.}}  (\ref{Legendre} $\Longrightarrow$ \ref{Legendre.normal}) Recalling the previous remarks, we only have to prove the sufficient condition [$\Longleftarrow$] of Theorem \ref{Legendre.normal}. Hence we consider an equation
\begin{equation}\label{eq1.legendre.normal.proof}
ax^2+by^2=z^2,
\end{equation}
with $a,b\in {\cal M}\models I\Delta_0$, $a,b$ square-free and positive, and suppose conditions $(\it Norm.1)\ a{\cal R}b$, $(\it Norm.2)\ b{\cal R}a$, and $(\it Norm.3)\ -{\frac{ab}{d^2}}{\cal R}d$ of Theorem \ref{Legendre.normal} hold, where $d=(a,b)$.

We now consider the equation
\begin{equation}\label{eq3.legendre.normal.proof}
Ax^2+By^2+Cz^2=0,
\end{equation}

where $A=\frac{a}{d}, B=\frac{b}{d}$ and $C=-d$. The coefficients $A,B,C$ are square-free and pairwise relatively prime. We have:
$$-AB=-\frac{ab}{d^2} \mbox{, and we know that } -{\frac{ab}{d^2}}{\cal R}d \mbox{ by  ({\it Norm.3})},$$

hence $-AB{\cal R} C$. Moreover,
$$-AC=-\frac{a}{d}(-d)=a, \mbox{ and we know that } a{\cal R} b \mbox{ by  ({\it Norm.3}), so } a{\cal R}\frac{b}{d},$$

that means $-AC{\cal R} B$. Finally,

$$-BC=-\frac{b}{d}(-d)=b \mbox{, and } b{\cal R} a \mbox{ by  ({\it Norm.2})},$$
which implies $b{\cal R}\frac{a}{d}$, that means $-BC{\cal R} A$.

So all conditions of Theorem \ref{Legendre} hold for the equation (\ref{eq3.legendre.normal.proof}). We can deduce that it has a non trivial solution $(x_0,y_0,z_0)$, i.e.
$$
\frac{a}{d}x_0^2+\frac{b}{d}y_0^2+(-d)z_0^2=0,
$$
which implies $ax_0^2+by_0^2=(dz_0)^2$.
Hence the triple $(x_0,y_0,dz_0)$ is a non-trivial solution of (\ref{eq1.legendre.normal.proof}).\\

(\ref{Legendre.normal} $\Longrightarrow$ \ref{Legendre}) As before, we only have to prove the sufficient condition [$\Longleftarrow$] of Theorem \ref{Legendre}.

Consider the equation
\begin{equation}\label{eq1.legendre.proof}
ax^2+by^2+cz^2=0,
\end{equation}
with $a,b,c\in {\cal M}\models I\Delta_0$, $a,b,c$ square-free, pairwise relatively prime and not all of the same sign. W.l.o.g. we can assume $a,b>0$ and $c<0$, and suppose conditions $(\it Leg.1)\ -ab{\cal R}c$, $(\it Leg.2)\ -bc{\cal R}a$, $(\it Leg.3)\ -ac{\cal R} b$
of Theorem \ref{Legendre} hold.

If we multiply both sides of (\ref{eq1.legendre.proof}) by $-c$ we obtain the equation
\begin{equation}\label{eq4.legendre.normal.proof}
Ax^2+By^2=Z^2,
\end{equation}
with $A=-ac, B=-bc$ and $Z=cz$. The coefficients $A$ and $B$ are square-free (since $a,b,c$ are pairwise relatively prime) and positive.

From ({\it Leg.3})  we have $A{\cal R} b$, and also $A{\cal R} c$ since $c|A$, and since $(b,c)=1$ we can apply Lemma \ref{lemmaCRT} to obtain that $A{\cal R} B$. In the same way we prove that $B{\cal R} A$.

Finally, notice that $(A,B)=(ac,bc)=c$, and
$$
-\frac{AB}{c^2}=1\frac{-(ac)(-bc)}{c^2}=-ab, \mbox{ and we know that } -ab{\cal R}c \mbox{ from (\it Leg.1)},
$$
hence we have  $-\frac{AB}{c^2}{\cal R} c$. We can conclude that all conditions of Theorem \ref{Legendre.normal} hold for the equation (\ref{eq4.legendre.normal.proof}), so it has a non trivial solution $(x_0,y_0,Z_0)$, for which $Z_0=cz_0$, for some $z_0\in\cal M$, and
$$
(-ac)x_0^2+(-bc)y_0^2=(cz_0)^2 \Longleftrightarrow c(-ax_0^2-by_0^2)=c^2z_0^2 \Longleftrightarrow ax_0^2+by_0^2=cz_0^2,
$$
hence $(x_0,y_0,z_0)$ is a non-trivial solution of (\ref{eq1.legendre.proof}).$_\Box$

\subsection{Proof of the theorem in $I\Delta_0+\Omega_1$}

In this section we are going to prove Legendre's theorem in the theory $I\Delta_0+\Omega_1$ by proving its equivalent normal form. We will pay close attention in adapting all the arguments of the corresponding proof in \cite{IR} to our theory.\\

A crucial point in the proof relies on the following property, which is well known for the integers. 

\begin{prop}\label{sumofsquares}
Let ${\cal M}\models I\Delta_0$ and $b\in {\cal M}$. If $-1$ is a square modulo $b$, then $b$ is the sum of two squares.
\end{prop}

On most texts on classic number theory this property is proved using tools such as the (full) pigeonhole principle, which is not available even in $I\Delta_0+\Omega_1$, or the relation about {\it Legendre's symbol} $\left(\frac{a}{p}\right)=a^{\frac{p-1}{2}}(mod\ p)$ (see \cite{HW}), which is not known to be valid in $I\Delta_0$. Here we show a different approach that fits to our context. We need some lemmas.

\begin{lem}\label{lemma1primesquare}
If $p\in{\cal M}\models I\Delta_0$ is a prime and $-1$ is a square modulo $p$ (i.e. $-1 {\mathcal R} p$), then there is $k\in{\cal M}$, with $k<p$, such that $kp=1+a^2$, for $a\leq \frac{p-1}{2}$.
\end{lem}

{\it{Proof.}}
Just notice that if $-1\equiv a^2(mod\ p)$, with $a\in\cal M$, $a\leq \frac{p-1}{2}$ be such, then there is $k\in \cal M$ such that
$$
kp=a^2+1 \leq \left(\frac{p-1}{2}\right)^2+1.
$$
It is easy shown that $\left(\frac{p-1}{2}\right)^2+1<p^2$, hence we have $kp<p^2$, with $k<p$, and $kp=1+a^2$. $_\Box$

\begin{lem}\label{lemma2primesquare}
Let $p\in{\cal M}\models I\Delta_0$ be a prime. If there is $k\in{\cal M}$, with $k<p$ such that $kp$ is the sum of two squares, then $p$ itself is the sum of two squares.
\end{lem}

{\it{Proof.}} It is clear that $2=1+1$ is a sum of two squares, so we can assume $p\not=2$.

The set $S(p)=\left\{ m\in{\cal M}: m< p \mbox{ and } mp \mbox{ is the sum of two squares }\right\}$, is $\Delta_0$-definable via the formula $\sigma(x) = \exists y<p\  \exists z<p\ (xp=y^2+z^2)$, and it is clearly bounded by $p$ and, by hypothesis, not empty. By $\Delta_0$-induction $S(p)$ has a minimum element, which we call $h$. So let $a,b<p$ be such that $hp=a^2+b^2$. Suppose that $h>1$, and we show that we can find a $h'\in S(p)$ with $h'<h$. We will distinguish two cases.

$\bullet$ {\it case $h$ is even}: then $a$ and $b$ must be both even or both odd.\\
If $a,b$ are both even we can write $hp=4\left(\left(\frac{a}{2}\right)^2+\left(\frac{b}{2}\right)^2\right)$,
hence $4$ divides $h$ and we have
$$\frac{h}{4}p=\left(\frac{a}{2}\right)^2+\left(\frac{b}{2}\right)^2.$$

If $a,b$ are both odd we can write
$$\frac{h}{2}p=\left(\frac{a-b}{2}\right)^2+\left(\frac{a+b}{2}\right)^2,$$
so in both cases we get $h'p$ as a sum of two squares, with $h'\leq h/2<h$, and this is a contradiction.

$\bullet$ {\it case $h$ is odd}: let $a\equiv \alpha(mod\ h)$ and $b\equiv\beta(mod\ h)$, with $\alpha,\beta<\frac{h}{2}$. Then
$$hp=a^2+b^2\equiv \alpha^2 +\beta^2\equiv 0(mod\ h),$$
hence there is $j\in{\cal M}$ such that
\begin{equation}\label{alphabeta}
 \alpha^2 +\beta^2=jh,
\end{equation}
and $j<h$ since $ \alpha^2 +\beta^2<(h/2)^2+(h/2)^2=h^2/2<h^2$.\\
Now from $a^2+b^2=kp$ and (\ref{alphabeta}) we obtain
\begin{equation}\label{alphabeta2}
\left(a^2+b^2\right)\left(\alpha^2+\beta^2\right)=jh^2p.
\end{equation}
Since
$$
\left(a^2+b^2\right)\left(\alpha^2+\beta^2\right)=\left(a\alpha+b\beta\right)^2+\left(a\beta-b\alpha\right)^2,
$$
and $a\alpha+b\beta\equiv\alpha^2+\beta^2\equiv 0(mod\ h)$, we have that $h$ divides $a\alpha+b\beta$ and, similarly, $h$ divides $a\beta-b\alpha$.\\
Henceforth, from (\ref{alphabeta2}) we obtain
$$
\left(\frac{a\alpha+b\beta}{h}\right)^2+\left(\frac{a\beta-b\alpha}{h}\right)^2=jp,
$$
so we have $jp$ as a sum of two squares and $j<h$, a contradiction.

It follows that $h=1$ and hence $p$ is a sum of two squares. $_\Box$

\medskip
It is now easy to deduce the following property.

\begin{prop}\label{primesumofsquares}
Let ${\cal M}\models I\Delta_0$ and let $p\in \cal M$ be a prime. If $-1$ is a square modulo $p$, then $p$ is the sum of two squares.
\end{prop}
{\it{Proof.}}
If $-1$ is a square modulo $p$, then by Lemma \ref{lemma1primesquare} there are $k,a\in \cal M$, $k,a<p$, such that $kp=1+a^2$, which is clearly a sum of two squares. The statement then follows straightforward from Lemma \ref{lemma2primesquare}.$_\Box$\\

\medskip

Proposition \ref{sumofsquares} is now easily obtained as follows.\\

{\it{Proof. (Proposition \ref{sumofsquares})}}
If $-1$ is a square modulo $b$, then $-1$ is a square modulo every prime $p$ dividing $b$. Then Proposition \ref{primesumofsquares} implies that every prime dividing $b$ is the sum of two squares.

It is easy to verify that the product of sums of two squares is still a sum of two squares (e.g. $(a^2+b^2)(c^2+d^2)=(ac+bd)^2+(ad-bc)^2$). So we obtain the result by iterating this argument to the product of all primes dividing $b$, which is of logarithmic length with respect to $b$ (roughly $log_2\ b$), where all the partial products are bounded by $b$ itself. $_\Box$

\medskip

We can now go through the proof of theorem \ref{Legendre.normal} in the theory $I\Delta_0+\Omega_1$. We will formalize by $\Delta_0$-induction a procedure to find a solution of the considered equation in models of $I\Delta_0+\Omega_1$. Here we restate the theorem.

\begin{theor}[Legendre normal form]\label{Legendre.normal.2}
Let ${\cal M}\models I\Delta_0+\Omega_1$ and let $a,b\in {\cal M}$ be square-free and positive.
Then the equation
\begin{equation}\label{eq.legendre.normal2}
ax^2+by^2=z^2
\end{equation}
has a non trivial solution in ${\cal M}$ if and only if
\begin{description}\vspace{-0.5em}
 \item{\rm (\it Norm.1)}\ $ a{\cal R}b$,
 \item{\rm (\it Norm.2)}\ $b{\cal R}a$,
 \item{\rm (\it Norm.3)}\ $-{\frac{ab}{d^2}}{\cal R}d,$ where $d=(a,b)$.
\end{description}
\end{theor}

{\it{Proof.}} We only have to prove that the properties ({\it Norm.1-3}) imply that the equation (\ref{eq.legendre.normal2}) has a non-trivial solution. Consider such an equation and suppose conditions ({\it Norm.1-3}) hold for $a,b\in{\cal M}$, $a,b$ square-free and positive. There are some trivial cases:

\begin{description}\vspace{-0.5em}
 \item{$\bullet$ Case $a=1$:}\ then (\ref{eq.legendre.normal2}) becomes $x^2+by^2=z^2$ and the triple $(1,0,1)$ is a non-trivial solution;
 \item{$\bullet$ Case $b=1$:}\ as before, being $(0,1,1)$ a non-trivial solution for $ax^2+y^2=z^2$;
 \item{$\bullet$ Case $a=b$:}\ then (\ref{eq.legendre.normal2}) becomes $b(x^2+y^2)=z^2$, and condition ({\it Norm.3}) becomes $-1{\cal R}b$. By Proposition \ref{sumofsquares} there are $r,s\in\cal M$ such that $b=r^2+s^s$, and so the triple $(r,s,b)$ represents a non-trivial solution of the equation.
\end{description}

Notice that in all these trivial cases the solution $(x_0,y_0,z_0)$ is such that $x_0,y_0,z_0\leq a$.

We can now consider $a,b>1$, $a\not= b$ and without loss of generality we suppose $b<a$. The argument is as follows: from the starting equation (\ref{eq.legendre.normal2}) we build another equation
\begin{equation}\label{eq.legendre.normal.A}
Ax^2+by^2=z^2,
\end{equation}
with $0<A<a$ and satisfying the appropriate conditions ({\it Norm.1-3}), such that if (\ref{eq.legendre.normal.A}) has a non-trivial solution, we obtain a non-trivial solution of (\ref{eq.legendre.normal2}) from it. By applying repeatedly this argument, and possibly switching the role of the coefficients $a$ and $b$ at some point, we eventually get to one of the trivial cases $a=1$, $b=1$ or $a=b$, that admit a non-trivial solution, and going backward from that we obtain a non-trivial solution to (\ref{eq.legendre.normal2}). We have to formalize this argument by $\Delta_0$ induction.

In the equation $ax^2+by^2=z^2$, ({\it Norm.2}) implies that there is $\beta\leq a/2$ such that $b\equiv\beta^2(mod\ a)$, hence there is $k\leq a$ such that $\beta^2-b=ka$. If we factor out the squares in $k$ we can write
\begin{equation}\label{newcoefficientA}
\beta^2-b=h^2Aa,
\end{equation}
with $A$ square-free. This is the coefficient we use in equation (\ref{eq.legendre.normal.A}) we are going to work with.

Froma \ref{newcoefficientA} it follows easily that $A>0$. We also get
$$
aA\leq aAh^2<\beta^2\leq\frac{a^2}{4},
$$
that means
\begin{equation}\label{Ainequality}
A<\frac{a}{4}.
\end{equation}
This inequality will turn out to be very important later.

Now let $d=(a,b)$, and $a=da_1,\ b=db_1$, with $(a_1,b_1)=1$. If a prime $p$ divides both $a_1$ and $d$, then $p^2$ divides $a$, and since $a$ is square-free we have $(a_1,d)=1$, and the same argument shows $(b_1,d)=1$.\\
From (\ref{newcoefficientA}) we get
$$
\beta^2=h^2Aa+b=h^2Aa_1d+b_1d=d(h^2Aa_1+b_1),
$$
hence $d$ divides $\beta^2$, and since $d$ is square-free, we have $d|\beta$. Llet $\beta=d\beta_1$, then $\beta^2=d^2\beta_1^2$ and we have $d^2\beta_1^2=d(h^2Aa_1+b_1)$, hence
\begin{equation}\label{dbeta1}
d\beta_1^2=h^2Aa_1+b_1.
\end{equation}
Now, if any prime $p$ divides both $d$ and $h$, from (\ref{dbeta1}) it follows that $p$ divides $b_1$, and hence $p$ divides both $b_1$ and $d$, a contradiction since they are relatively prime. So necessarily $(d,h)=1$. From (\ref{dbeta1}) we obtain
$$
h^2Aa_1\equiv -b_1(mod\ d)\ \mbox{ which implies }\ h^2Aa_1^2\equiv -b_1a_1(mod\ d),
$$
and since $(a_1,d)=(h,d)=1$,  both $h$ and $a_1$ are invertible modulo $d$. So we have
$$
A\equiv -a_1b_1\left(h^{-1}a_1^{-1}\right)^2(mod\ d).
$$
Now, since $-a_1b_1=-\frac{ab}{d^2}$, condition ({\it Norm.3}) tells us that $-a_1b_1$ is a square modulo $d$, and so we get
\begin{equation}\label{Asquare1}
A{\mathcal R}d.
\end{equation}
Notice that if there is a prime $p$ which divides both $h$ and $b$, then from (\ref{newcoefficientA}) we get that $p$ divides $\beta$ an then $p^2$ divides $b$, and this is a contradiction since $b$ is square-free. Hence $(h,b)=1$. Since $b_1$ divides $b$, also $(h,b_1)=1$. Then from (\ref{newcoefficientA}) the following implication holds
$$
h^2Aa\equiv\beta^2(mod\ b_1) \Rightarrow A\equiv\beta^2\left(h^{-1}\right)^2a^{-1}(mod\ b_1).
$$
From ({\it Norm.1}) states that $a{\cal R}b$ it follows that $a{\cal R}b_1$, and so
\begin{equation}\label{Asquare2}
A{\mathcal R}b_1.
\end{equation}
Now from (\ref{Asquare1}), (\ref{Asquare2}) and Lemma \ref{lemmaCRT} we get
\begin{equation}\label{norm.1A}
A{\mathcal R}b.
\end{equation}
Moreover, from (\ref{newcoefficientA}) we know that $b\equiv\beta^2(mod\ A)$, that means
\begin{equation}\label{norm.2A}
b{\mathcal R}A.
\end{equation}

If we put $r=(A,b)$, it is left to show that $-\frac{Ab}{r^2}{\cal R}r$.\\
Let $A=A_2r, b=b_2r$, with $(A_2,b_2)=1$. We have $(r,A_2)=(r,b_2)=1$ since $A,b$ are square-free. From (\ref{newcoefficientA}) we obtain
\begin{equation}\label{betar}
\beta^2=b_2r+h^2A_2ra=r(b_2+h^2A_2a).
\end{equation}
So $r$ divides $\beta^2$, and since $r$ is square-free, we have $r|\beta$. Let $\beta=\beta_2r$, from (\ref{betar}) we get the following implications
\begin{equation}\label{A2square1}
r\beta_2^2=b_2+h^2A_2a \Rightarrow h^2A_2a\equiv -b_2(mod\ r)\Rightarrow -A_2b_2h^2a\equiv b_2^2(mod\ r).
\end{equation}
Now using the same arguments as before we can show that $(a,r)=(h,r)=1$, so both $a,h$ are invertible modulo $r$, and we obtain
$$
-A_2b_2\equiv b_2^2\left(h^{-1}\right)^2a^{-1}(mod\ r).
$$
Recalling that $a{\cal R}b$ and $r|b$, we have that $a{\cal R}r$, and since $-A_2b_2=-\frac{Ab}{r^2}$, the previous congruences imply that
\begin{equation}\label{norm.3A}
-\frac{Ab}{r^2}{\mathcal R}r.
\end{equation}

We have then obtained the equation (\ref{eq.legendre.normal.A}) $Ax^2+by^2=z^2$, with $0<A<a/4$, $A$ square-free and (by (\ref{norm.1A}), (\ref{norm.2A}) and (\ref{norm.3A})) satisfying the conditions
\begin{description}\vspace{-0.5em}
 \item{\rm (\it Norm.1)}\ $ A{\cal R}b$,
 \item{\rm (\it Norm.2)}\ $b{\cal R}A$,
 \item{\rm (\it Norm.3)}\ $-{\frac{Ab}{r^2}}{\cal R}r,$ where $r=(A,b)$.
\end{description}

The single reduction we have made can easily be formalized in $I\Delta_0$, since it is based only on congruences and the quantifiers are obviously bounded by the initial coefficients $a$ and $b$.

Now suppose $(x_0,y_0,z_0)$ is a non-trivial solution of (\ref{eq.legendre.normal.A}), hence
\begin{equation}\label{solutionA}
Ax_0^2=z_0^2-by_0^2
\end{equation}
By multiplying (\ref{solutionA}) by (\ref{newcoefficientA}) we have
$$
A^2x_0^2h^2a=(z_0^2-by_0^2)(\beta^2-b)=z_0^2\beta^2-bz_0^2-by_0^2\beta^2+b^2y_0^2;
$$
if we now add and subtract the quantity $2z_0\beta by_0$ we have
$$
A^2x_0^2h^2a=(z_0^2\beta^2+b^2y_0^2+2z_0\beta by_0)-b(z_0^2+y_0^2\beta^2+2z_0\beta by_0)=(z_0\beta+by_0)^2-b(z_0+y_0\beta)^2,
$$
and so
$$
a(Ax_0h)^2 + b(z_0+y_0\beta)^2 =(z_0\beta+by_0)^2,
$$
which states that the triple
$$
(Ax_0h,\ z_0+y_0\beta,\ z_0\beta+by_0)
$$
is a non-trivial solution of equation (\ref{eq.legendre.normal2}).

We now need to estimate the growth rate of the solution of equation (\ref{eq.legendre.normal2}) in terms of that of (\ref{eq.legendre.normal.A}). First of all, notice that if $Ax_0^2+by_0^2=z_0^2$, then clearly $x_0,y_0\leq z_0$ (remember that both $A$ and $b$ are positive). Then, as already showed, the components of the solution of (\ref{eq.legendre.normal2}) are
\begin{description}\vspace{-0.5em}
 \item $\bullet\ Ax_0h\leq x_o\frac{a}{4} \mbox{ (since } Ah<\frac{a}{4}{\mbox)}$
 \item $\bullet\ z_0+y_0\beta\leq  z_0+z_0\beta=z_0\left(1+\frac{a}{2}\right) \mbox{ (since } \beta\leq\frac{a}{2}{\mbox)}$
 \item $\bullet\ z_0\beta+by_0\leq z_0\frac{a}{2}+az_0=z_0\left(\frac{3}{2}a\right)$
\end{description}
and since we are assuming $a>1$ we can conclude that all the components of the new solution are $\leq z_0\left(\frac{3}{2}a\right)$.\\

We now have to iterate this procedure and formalize it in $I\Delta_0+\Omega_1$. We start with the given equation
\begin{equation}\label{E0}
E_0:\ ax^2+by^2=z^2,
\end{equation}
where $a,b$ are square-free, $a>b$ and
$$
a{\cal R}b,\ b{\cal R}a,\ -\frac{ab}{d^2}{\cal R}d,\mbox{ where } d=(a,b),
$$
and we build a sequence of equations, for $i>0$
\begin{equation}\label{Ei}
E_i:\ A_ix^2+B_iy^2=z^2,
\end{equation}
where every $A_i,B_i$ are defined by recursion as follows (where $A_0=a, B_0=b$)
\begin{equation}\label{betai}
\bullet\ \beta_i^2-B_i=h_i^2A_iA_{i-1},\ B_i=B_{i-1}
\end{equation}
with $A_ih_i<\frac{A_{i-1}}{4},\ \beta_i\leq \frac{A_{i-1}}{2}$, if at step $i-1$ we have $A_{i-1}>B_{i-1}$, or
\begin{equation}\label{alphai}
\bullet\ \alpha_i^2-A_i=h_i^2B_iB_{i-1},\ A_i=A_{i-1}
\end{equation}
with $B_ih_i<\frac{B_{i-1}}{4},\ \alpha_i\leq \frac{B_{i-1}}{2}$, if at step $i-1$ we have $A_{i-1}<B_{i-1}$.

For every equation $E_i$ the following congruence conditions hold.
$$
A_i{\cal R}B_i,\ B_i{\cal R}A_i,\ -\frac{A_iB_i}{r_i^2}{\cal R}r_i,\mbox{ where } r_i=(A_i,B_i).
$$
Moreover, if $(x_i,y_i,z_i)$ is a non-trivial solution of equation $E_i$, then a non-trivial solution of $E_{i-i}$ is either
\begin{equation}\label{solutionbetai}
\bullet\ (A_ix_ih_i,\ z_i+y_i\beta_i,\ z_i\beta_i+B_{i-1}y_i)
\end{equation}
or
\begin{equation}\label{solutionalphai}
\bullet\ (z_i+y_i\alpha_i,\ B_ix_ih_i,\ z_i\alpha_i+A_{i-1}y_i)
\end{equation}
according to $A_{i-1}>B_{i-1}$ or  $A_{i-1}<B_{i-1}$, respectively.

We remark that:
\begin{description}\vspace{-0.5em}
	\item {\it{(i)}} the growth factor from a solution of the equation $E_i$ to that of $E_{i-1}$ is always bounded by $\frac{3}{2}A_{i-1}\leq\frac{3}{2}a$ when $A_{i-1}>B_{i-1}$, and by $\frac{3}{2}B_{i-1}\leq\frac{3}{2}b$ when $A_{i-1}<B_{i-1}$;
	\item {\it{(ii)}} when "descending" through the sequence, the coefficients of equation $E_i$ and those of equation $E_{i-1}$ are related as follows:
	$$A_i<\frac{A_{i-1}}{4}\leq\frac{a}{4} \mbox{ or } B_i<\frac{B_{i-1}}{4}\leq\frac{b}{4}.$$ Hence the length of the sequence of $E_i$'s is at most $log_{_4}a+log_{_4}b$.
	\item {\it{(iii)}} the sequence of equations will eventually stop with one of the trivial cases where one of the coefficients is 1 or they are equal. In these cases non trivial solutions exist, namely $(1,0,1)$ or $(0,1,1)$ or $(r_i,s_i,B_i)$, where $B_i=r_i^2+s_i^2$.
\end{description}

Let $l$ be the length of the sequence. For the final solution of equation $E_l$ we can clearly state that $x_l,y_l,z_l\leq b$ (where $b$ is the coefficient of the initial equation (\ref{eq.legendre.normal2})). From remarks {\it{(i)}} and  {\it{(ii)}} we can derive that all the components of the non-trivial solution $(x_0,y_0,z_0)$ of  $E_0$, and so of  (\ref{eq.legendre.normal2}), are bounded as follows
\begin{equation}\label{bound}
x_0,y_0,z_0\leq b\left(\frac{3}{2}a\right)^{log_{_4}a}\left(\frac{3}{2}b\right)^{log_{_4}b}.
\end{equation}

It is only left to formalize the recursion we have constructed in $I\Delta_0+\Omega_1$.\\
For the sake of clarity we will recall here all the $\Delta_0$-formulas we need:
\begin{description}\vspace{-0.5em}
	\item $\bullet$ $x$ divides $y$:  $\delta(x,y)=\exists z\leq y\ (xz=y)$
	\item $\bullet$ $x$ is a prime: $Pr(x)=\forall y\leq x\ (\delta(y,x)\rightarrow (y=1\lor y=x))$
	\item $\bullet$ $z$ is the $g.c.d.$ of $x$ and $y$: $$\gamma(x,y,z)=\delta(z,x)\land\delta(z,y)\land \forall t\leq x\ (\delta(t,x)\land\delta(t,y))\rightarrow \delta(t,z)$$
	\item $\bullet$ $x$ is square-free: $\sigma(x)=\forall y\leq x\ (Pr(y)\rightarrow \lnot\delta(y^2,x))$
	\item $\bullet$ $x$ is a square modulo $y$: $\rho(x,y)=\exists z\leq x \land \exists r\leq y/2\ (x=r^2+zy)$.
\end{description}

We can now express all conditions of the theorem with a $\Delta_0$-formula:
\begin{equation}\label{conditions}
\begin{split}
\Theta(a,b)= &\  0\leq a\land 0\leq b\land \sigma(a)\land \sigma(b)\land \rho(a,b)\land \rho(b,a) \land\\
 & \land \forall d\leq a\  \left(\gamma(a,b,d)\rightarrow \rho\left(-\frac{ab}{d^2},d\right)\right).
\end{split}
\end{equation}

We now make induction on the formula
\begin{equation*}
\Lambda(t)= \forall a\leq t\ \forall b\leq t\ \left(ab\leq t \land b\leq a\land \Theta(a,b)\right)\longrightarrow 
\end{equation*}
\begin{equation}\label{induction}
\exists x,y,z\leq b\left(\frac{3}{2}a\right)^{log_{_4}a}\left(\frac{3}{2}b\right)^{log_{_4}b}\ \lnot(x=0\land y=0\land z=0)\land (ax^2+by^2=z^2),
\end{equation}
which is a $\Delta_0$ formula that uses boundaries which are allowed by the axiom $\Omega_1$.

For $t=1$ the formula is true since in this case $a=b=1$ and the equation $x^2+y^2=z^2$ has non-trivial solutions, for example $(1,0,1)$, which clearly satisfy $x,y,z\leq 1\cdot\left(\frac{3}{2}\cdot1\right)^{log_{_4}1}\left(\frac{3}{2}\cdot1\right)^{log_{_4}1}=1$.

Now suppose ${\cal M}\models \Lambda(t)$, with $t\in\cal M$, $t>1$, and consider $t'=t+1$\\
Let $a,b\in \cal M$, $a,b\leq t'$;  and $ab=t'$ (if $ab<t'\Rightarrow ab\leq t$ and we already know $\Lambda(t)$ is true).\\
W.l.o.g. we can assume $b<a$, and hence apply the first step of the reduction and obtain the equation $Ax^2+by^2=z^2$, with $A<a/4$ and ${\cal M}\models \Theta(A,b)$.

Now $Ab<t'$, hence $Ab\leq t$, so by inductive hypothesis (${\cal M}\models \Lambda(t)$), this equation admits a non-trivial solution $(x_1,y_1,z_1)$ in $\cal M$, such that
$$
x_1,y_1,z_1\leq b\left(\frac{3}{2}A\right)^{log_{_4}A}\left(\frac{3}{2}b\right)^{log_{_4}b}.
$$

From this we showed how we can get a non-trivial solution $(x_0,y_0,z_0)$ of the equation $ax^2+by^2=z^2$, and by the previous observations we made we can state that
$$
x_0,y_0,z_0\leq b\left(\frac{3}{2}A\right)^{log_{_4}A}\left(\frac{3}{2}b\right)^{log_{_4}b}\left(\frac{3}{2}a\right)\leq b\left(\frac{3}{2}a\right)^{log_{_4}A+1}\left(\frac{3}{2}b\right)^{log_{_4}b}.
$$
Since $A<a/4$, we have $log_{_4}A\leq log_{_4}a-1$, and we can conclude that
$$
x_0,y_0,z_0\leq b\left(\frac{3}{2}a\right)^{log_{_4}a}\left(\frac{3}{2}b\right)^{log_{_4}b}.
$$
Hence we have that ${\cal M}\models\Lambda(t')$, and this concludes the proof. $_\Box$

\medskip

{\bf Concluding remarks:}

We have adapted a proof of Legendre's theorem suggested in \cite{IR}. The proof we exhibit provides a bound for the solution of the initial equation which is in terms of $b\left(\frac{3}{2}a\right)^{log_{_4}a}\left(\frac{3}{2}b\right)^{log_{_4}b}$. This is the reason why we can state Legendre's theorem only in $I\Delta_0+\Omega_1$, even though all previous properties and statements are valid even in $I\Delta_0$.
Cassels in \cite{Cassels} obtained a linear bound of the solution in terms of the initial coefficients. Unfortunately the proof uses tools of geometry of numbers which seem to rely on the (full) pigeonhole principle, which is not known to be provable in $I\Delta_0+\Omega_1$. Hence a possible further development in this subject could be to search for an alternative proof with no use of $PHP$, in order to obtain Cassels' result in $I\Delta_0+\Omega_1$ or even in $I\Delta_0$.


\begin{thebibliography}{x}
\addcontentsline{toc}{chapter}{Bibliography}

\bibitem[BI]{Ber-Intr}
Berarducci A. and Intrigila B., {\it Combinatorial principles in elementary number theory}, Annals of Pure and Applied Logic, vol. 55 (1991), pp. 35-50. 

\bibitem[C]{Cassels}
Cassels J.W.S., {\it Rational Quadratic Forms}, Academic Press, 1978.

\bibitem[D]{D'Aquino}
D'Aquino P., {\it Pell equations and exponentiation in fragments of arithmetic}, Annals of Pure and Applied Logic, vol. 77 (1996), pp.1-34.

\bibitem[D2]{D'Aquino2}
D'Aquino P., {\it Weak fragments of Peano Arithmetic}, in The Notre Dame Lectures, edited by P. Cholak, Association for Symbolic Logic, Lecture Notes in Logic, 18 (2005), pp. 149-185

\bibitem[D3]{D'Aquino3}
D'Aquino P., {\it Local behaviour of Chebyshev teorem in models of $I\Delta_0$}, Journal of Symbolic Logic 57 (1) (1992), pp.12-27

\bibitem[DM]{D-M}
D'Aquino P. and A. Macintyre, {\it Non standard finite fields over $I\Delta_0+\Omega_1$}, in Israel Journal of Mathematics 117, (2000), pp. 311-333.

\bibitem[DM2]{DM2}
D'Aquino P. and A. Macintyre, {\it Primes in Models of $I\Delta_0+\Omega_1$: Density in Henselizations}, in New Studies in Weak Arithmetics, edited by P. Cegielski, C.Cornaros, C. Dimitracopoulos, CSLI Publications, Lecture Notes Number 211, pp. 85-91.

\bibitem[DM3]{DM3}
D'Aquino P. and A. Macintyre, {\it Quotient Fields of a Model of  $I\Delta_0+\Omega_1$}, Mathematical Logic Quarterly 47 (2001) 3, pp. 305-314.

\bibitem[FdG]{F-dG}
Franciosi S., de Giovanni F., {\it Elementi di Algebra}, Aracne ed., 1995.

\bibitem[GD]{Gaifman-Dimitracopulos}
Gaifman H., Dimitracopulos C., {\it Fragments of Peano's arithmetic and the MRDP theorem}, Logic and Algorithmic (Zurich, 1980), Univ. Genève, Geneva, 1982, pp. 187-206.

\bibitem[HP]{HP}
Hajek P. and Pudlak P., {\it Metamathematics of first-order arithmetic}, Springer-Verlag, Berlin, 1998, second printing.

\bibitem[HW]{HW}
Hardy G.H. and Wright E.M., {\it An Introduction to the Theory of Numbers}, (5th ed.) Oxford University Press, 1980.

\bibitem[IR]{IR}
Ireland K., Rosen M., {\it A Classical Introduction to Modern Number Theory}, second edition, Springer, 1990.

\bibitem[K]{Kaye}
Kaye R., {\it Models of Peano Arithmetic}, Oxford University Press, Oxford, 1991.

\bibitem[L]{Lagarias}
Lagarias J.C., {\it On the computational complexity of determining the solvability or unsolvability of the equation $X^2-DY^2=-1$}, in Transactions of American Mathematical
Society vol. 260, N. 2 (1980), pp. 485-508.

\bibitem[MM]{MM}
Macintyre A., Marker D., {\it Primes and their residue rings in models of open induction}, Annals of Pure and Applied Logic, vol. 43 (1989), no. 1, pp 57-77.

\bibitem[MA]{A-M}
Manders K. L., Adleman L., {\it NP-complete decision problems for binary quadratics}, Journal of Computer and System Sciences, vol. 16 (1978), no. 2, pp. 168-184.

\bibitem[Mar]{Marker}
Marker D., {\it Model Theory: an Introduction}, Springer-Verlag, 2002

\bibitem[Ot]{Ot}
Otero M., {\it Models of open induction}, Ph.D. Thesis, Oxford University, 1991.

\bibitem[Pa]{Parikh}
Parikh R., {\it Existence and feasibility in arithmetic}, Journal of Symbolic Logic, vol. 36 (1976), no. 3, pp 494-508. 

\bibitem[PWW]{Paris-Wilkie-Woods}
Paris J., Wilkie A. and Woods A., {\it Provability of the Pigeonhole Principle and the existence of infinitely many primes}, Journal of Symbolic Logic 53, no. 4,
(1988), pp 1235-1244. 

\bibitem[R]{Rose}
Rose H.E., {\it A Course in Number Theory}, 2nd ed. Oxford University Press, 1995.

\bibitem[S]{Shep}
Shepherdson C.J., {\it A non-standard model for a free variable fragment of number theory}, Bull. Acad. Polon. Sci. Sér. Sci. Math. Astronom. Phys., vol. 12 (1964), pp. 79-86.

\bibitem[W]{Wilkie}
Wilkie A., {\it Applications of complexity theory to $\Sigma_0$-definability problems in arithmetic}, in Pacholski et al. eds., Model theory, Algebra and Arithmetic, Proc. Karpacz, Poland 1979. Lecture Notes in Mathematics vol. 834, Springer 1980, pp. pp.363-369.

\bibitem[WP]{Wilkie-Paris}
Wilkie A. and Paris J., {\it On the scheme of induction for bounded arithmetic formulas}, Annals of Pure and Applied Logic 35  (1987), pp. 261-302

\bibitem[Wo]{Woods}
Woods A., {\it Some problems in logic and number theory and their connections}, Ph.D. thesis, Manchester University, 1981.

\end{thebibliography}
\end{document}